\scrollmode
\documentstyle[11pt,amssymb]{article}
\def\a{{\alpha}}
\def\b{{\beta}}
\def\eps{{\varepsilon}}
\def\om{{\omega}}
\def\p{{\varphi}}
\def\l{{\lambda}}
\def\d{\partial}
\def\D{{\nabla}}
\def\bib#1{\bibitem[#1]{#1}}
\def\<{\langle}
\def\>{\rangle}
\def\id{\mbox{\rm id}}
\def\ad{\mathop{\mbox{\rm ad}}}
\def\diag{\mathop{\mbox{\rm diag}}}
\def\R{{\mathbb R}}
\def\N{{\mathbb N}}
\def\Z{{\mathbb Z}}
\def\O{{\Omega}}
\def\tf#1#2{{\textstyle\frac{#1}{#2}}}
\def\re#1{(\ref{#1})}
\def\eas{\begin{eqnarray*}}
\def\eeas{\end{eqnarray*}}
\def\qed{\hfill$\Box$\\[-1mm]}
\def\proof{\noindent{\bf Proof. }}
\def\nn{\nonumber}
\def\YM{Y\!M}
\def\so{{\mathfrak so}}
\def\xq{{\bar{x}}}
\newtheorem{definition}{Definition}

\newtheorem{proposition}{Proposition}
\newtheorem{lemma}{Lemma}
\newtheorem{remark}{Remark}
\newtheorem{theorem}{Theorem}

\begin{document}
\sloppy
\def\today{}
\title{Yang-Mills connections of cohomogeneity one\\
on $SO(n)$-bundles over Euclidean spheres}
\author{Andreas Gastel}
\maketitle

\section{Introduction}

For many geometric variational problems or p.d.e., there is a construction of
spherical solutions via {\em joins of spheres. }Very roughly, these 
constructions use the fact that $S^{m+n+1}$ is the {\em join\/} of $S^m$
and $S^n$, which means a set made up from curves, each one of which connects
one point in $S^m$ with one in $S^n$. These curves allow a common 
parametrization over $[0,\frac{\pi}2]$, say, and depending on this parameter
$t$, one can try to construct all kinds of geometric objects on
$S^{m+n+1}$ from homogeneous objects of the same type on $S^m$ and $S^n$.
Homogeneity of the latter helps reducing the partial differential equations,
which usually describe such objects, to {\em ordinary\/} differential 
equations.
The symmetries described here often, but not always, correspond to
some $SO(m+1)\times SO(n+1)$-invariance or -equivariance of the objects
being constructed.
This family of constructions has lead to examples of 
\begin{itemize}
\item harmonic maps between spheres \cite{Sm}, \cite{Di}, \cite{PR} (and also
some variants like $p$-harmonic maps \cite{Fa} and biharmonic maps \cite{GZ});
\item constant mean curvature hypersurfaces in $\R^{n+1}$ \cite{Hs1};
\item non-equatorial minimal embeddings of $S^n$ in $S^{n+1}$
  \cite{Hs2}, \cite{Hs3};
\item Einstein metrics on spheres \cite{Bo}.
\end{itemize}
A common feature of these constructions is that they all reduce the original
p.d.e. to a (system of) o.d.e.\ with {\em singular boundary values. }They 
tend to 
work best in dimensions which are slightly above the ``critical dimension'' 
of the respective equation. An excellent source presenting the first three of 
the examples in a unified way is the book \cite{ER2}.

The aim of this paper is to establish a similar construction for Yang-Mills
connections; more precisely for Yang-Mills $SO(n)$-connections over some
$S^m$, $m\ge5$. The join construction for such connections will exhibit
features similar to the ones listed above. Note that, due to the supercritical
dimension $m\ge5$ and to the fact that we cannot work in a ``Hermitian
Yang-Mills'' setting, there is currently no way to prove existence of such
connections by variational methods. This is of course closely related to the
lack of good gauges for connections in these dimensions. We make up for this
by choosing a suitable equivariant ansatz which already is in a ``good''
gauge. 

The methods we use are close to the methods for harmonic maps invented
by Smith, and the conditions for solvability are very reminiscent of the
``damping conditions'' known from harmonic map theory. 
The most notable difference is that in our case the equivariant
ansatz does not reduce the problem to a single o.d.e., but to a system of two
o.d.e. This fact adds a little bit of the flavor of B\"ohm's construction of
Einstein metrics to our considerations.

What we are going to construct are Yang-Mills connections of cohomogeneity 
one. It
should be noted that such connections over manifolds of dimension four have
been studied extensively by Urakawa \cite{Ur}, who also provides a very
general reduction setting for o.d.e.\ in Yang-Mills theory. The typical
degree three nonlinearities for the o.d.e.s found there and here probably
appeared first in Parker's construction of non-minimizing Yang-Mills
fields \cite{Par}.
Recently, Park and Urakawa \cite{PU} have also studied completely homogeneous 
Yang-Mills connections, which in a special case we will also have to do in 
this paper.

The paper is organized as follows. In {\bf Section 2}, we give a rather general
short introduction to equivariant Yang-Mills connections. {\bf Section 3} is
devoted to a rather more special case of homogeneous pull-back bundles
of $TS^n$ under mappings $S^m\to S^n$. We will need the so-called
``Yang-Mills eigenmaps'' obtained from these consideration as the homogeneous
``building blocks'' for our join construction.

In {\bf Section 4}, we observe that there is only a very restricted class
of joins of {\em vector bundles\/} which are again smooth vector bundles.
This justifies our reduction ansatz in {\bf Section 5}, which otherwise would
look a bit special at first glimpse. In this section, the reduction of the 
Yang-Mills
equation to a system of two o.d.e.\ (equipped with singular boundary data)
is performed.

The solvability of the singular o.d.e.\ boundary value problem thus 
obtained is 
discussed in some detail in {\bf Section 6}. We get sufficient conditions
for that, which (comparing to the harmonic map case where this is known)
we expect to be also necessary. Finally, in {\bf Section 7}, we apply the
existence theorem to find nontrivial examples of smooth cohomogeneity one
Yang-Mills connections over spheres. 
Among the examples we construct are
\begin{itemize}
\item one Yang-Mills connection on each of the countably many principal
  $SO(6)$-bundles over $S^6$,
\item countably many Yang Mills connections on $TS^n$ for 
  $n\in\{5,\ldots,9\}$.
\end{itemize}
(Coincidentally, B\"ohm's join construction for Einstein metrics \cite{Bo}
produces nonhomogeneous Einstein metrics on $S^n$ for exactly the same range
of dimensions.)

{\bf Acknowledgment:} This paper was finished while the author was visiting
SFB Transregio 71 at Freiburg. He would like to thank for hospitality and
support.

\section{Equivariant connections}

Let $M$ be a compact Riemannian manifold,
$\pi:E\to M$ be a $G$-vectorbundle of rank $n$ for some compact Lie group
$G\subseteq O(n)$; the latter means that we view $E$ equipped with a bundle
metric.

We assume that another compact Lie group $K$ is acting on both $E$ and $M$ by
isometries; we denote the action of $k\in K$ on $M$ simply by $k:M\to M$,
while on $E$ we denote it by $\l_k:E\to E$. We assume that the $K$-actions
are compatible with the projection, which means
\[
  \pi(\l_kv)=k\pi(v)
\]
for all $k\in K$, $v\in E$. By $\O^0(E)$
we denote the set of smooth sections of $E$, and by $\O^\ell(E)$, 
$\ell\in\N\cup\{0,\infty\}$, the sections
of $E\otimes\wedge^\ell T^*M$, i.e.\ the corresponding section-valued 
$\ell$-forms. 

The $K$-actions introduced above
induce a natural $K$-action on $\O^0(E)$, with $\tau_k:\O^0(E)
\to\O^0(E)$ given for $k\in K$ by
\[
  (\tau_k Y)(x):=\lambda_kY(k^{-1}x).
\]
A connection $D:\O^0(E)\to\O^1(E)$ is called {\em $K$-equivariant\/} if
\[
  D_{k_*u}(\tau_k Y)=\tau_k(D_uY)
\]
holds for all $k\in K$, $u\in\O^0(TM)$, and $Y\in\O^0(E)$. Here $k_*$ means 
the derivative of $k:M\to M$. 

Let us fix a $K$-equivariant reference connection $\D$ of $E$. Then every
$G$-connection of $E$ is of the form $D=\D+A$ for some $A\in\O^1(\ad P)$,
where $P$ is the principal fiber bundle associated with $E$.
We want to describe what equivariance of $D$ (and $\D$) means for $A$.
For $k,u,Y$ as above, we have
\eas
  \D_uY(x)+A_u(x)Y(x)&=&D_uY(x)\\
  &=&\tau_k^{-1}(D_{k_*u}(\tau_kY))(x)\\
  &=&\tau_k^{-1}(\D_{k_*u}(\tau_kY))(x)+\tau_k^{-1}(A_{k_*u}\tau_kY)(x)\\
  &=&\D_uY(x)+\tau_k^{-1}(A_{k_*u}\tau_kY)(x)\\
  &=&\D_uY(x)+\l_k^{-1}(A_{k_*u}\tau_kY)(kx)\\
  &=&\D_uY(x)+\l_k^{-1}A_{k_*u}(kx)\l_kY(x),
\eeas
from which we read off that
\[
  A_u(x)=\l_k^{-1}A_{k_*u}(kx)\l_k
\]
for all $x\in M$, $u\in E_x$, and $k\in K$. 
Similarly, we find the correct transformation of the curvature $F=F_A$ 
of $D$:
\eas
  F_{uv}(x)Y(x)&=&(D_uD_v-D_vD_u)Y(x)\\
  &=&\tau_k^{-1}(D_{k_*u}D_{k_*v}-D_{k_*v}D_{k_*u})(\tau_kY)(x)\\
  &=&\tau_k^{-1}(F_{k_*u,k_*v}\tau_kY)(x)\\
  &=&\l_k^{-1}F_{k_*u,k_*v}(kx)\l_kY(x),
\eeas
and hence
\[
  F_{uv}(x)=\l_k^{-1}F_{k_*u,k_*v}(kx)\l_k
\]
for all $x\in M$, $u,v\in E_x$, $k\in K$.

A connection $D_A=\D+A$
is called a {\em Yang-Mills connection\/}, if it is a critical
point of the Yang-Mills functional
\[
  \YM(A)=\frac12\int_M|F_A|^2\,dx.
\]
A connection is Yang-Mills if and only if
\[
  D_A^*F_A=0,
\]
which for smooth $A$ is equivalent to the weak formulation
\[
  \int_M\<F_A,D_A\p\>\,dx=0\qquad
  \mbox{ for all }\p\in\O^1(\ad P).
\]
A first important observation about equivariant Yang-Mills maps is an
instance of Palais' so-called {\em principle of symmetric criticality}, cf.\
\cite{Pal} for the general philosophy.

\begin{proposition}[symmetric criticality]\label{psc}
A smooth $K$-equivariant connection $D_A$ on $E$ is already Yang-Mills
if it is only critical with respect to equivariant variations, i.e.\ if
the first variation
\[
  \int_M\<F_A,D_A\p\>\,dx=0
\]
vanishes for those $\p\in\O^1(\ad P)$ satisfying
\[
  \p_u(x)=\l_k^{-1}\p_{k_*u}(kx)\l_k
\]
for all $x\in M$, $u\in E_x$, $k\in K$.
\end{proposition}

\proof
We abbreviate the right-hand side of the last equation by $(k^*\p)(x)$, and
similarly for $F$. Let $\p\in\O^1(\ad P)$ be {\em any\/} form, not necessarily
equivariant. Denoting the Haar measure of $K$ by $H_K$, and using the fact
that all $K$-actions are isometric and commute with $D$, we calculate
\eas
  \int_M\<F_A,D_A\p\>\,dx
  &=&\int_M\int_K\<(k^{-1})^*F_A,D_A\p\>\,dH_K\,dx\\
  &=&\int_M\int_K\<F_A,k^*(D_A\p)\>\,dH_K\,dx\\
  &=&\int_M\Big\<F_A,D_A\int_Kk^*\p\,dH_K\Big\>\,dx\\
  &=&0,
\eeas
where the first ``$=$'' holds because $F_A$ is $K$-equivariant, and the 
last one because so is $\int_Kk^*\p\,dH_K$. This proves that $D_A$
is Yang-Mills.\qed

\section{Homogeneous connections over $S^m$}

We start with some notation.
For $a,b\in\R^n$, we denote by $a\otimes b:\R^n\to\R^n$ the linear mapping
given by
\[
  (a\otimes b)(v):=\<a,v\>b,
\]
represented by the matrix
\[
  (a\otimes b)_{ij}=a_jb_i.
\]
If $M$ is a skew-symmetric $n\times n$-matrix, we have
\[
  (a\otimes b)_{ij}M_{jk}=a_jM_{jk}b_i=-M_{kj}a_jb_i
\]
and hence
\[
  (a\otimes b)M=-(Ma)\otimes b.
\]
Similarly,
\[
  M_{ij}(a\otimes b)_{jk}=M_{ij}a_kb_j,
\]
which means
\[
  M(a\otimes b)=a\otimes(Mb).
\]

We want to consider equivariant bundles over $S^m$ on the pull-back bundle
$E=h^*TS^n$ of some smooth map $h:S^m\to S^n$. We describe the bundle globally
by identifying the fiber over $x\in S^m$ with $T_{h(x)}S^n$, i.e. we
identify its total space with
\[
  E=\{(x,y)\in S^m\times\R^{n+1}:\<h(x),y\>=0\}.
\]
The total space
of the corresponding principal fiber bundle $P$ can be identified with
\[
  \{(x,M)\in S^m\times SO(n+1):Mh(x)=h(x)\},
\]
where the fiber over $x$ is the isotropy subgroup of $h(x)$ in $O(n+1)$ acting
on $\R^{n+1}\supset S^n$ in the standard way. Every connection on $P$
(or equivalently on $h^*TS^n$)
is of the form
\[
  D=\D+A,
\]
where here $A$ is a section in the adjoint vector bundle $\ad P$ 
with total space
\[
  \{(x,A)\in S^m\times\so(n+1):Ah(x)=0\},
\]
and $\D$ is the pull-back of the Levi-Civita connection on $S^n$. The 
latter means
\[
  \D_uY=\d_u Y-\<h,\d_u Y\>h=\d_u Y+\<\d_u h,Y\>h
\]
for all sections $u$ of $TS^m$ and $Y$ of $h^*TS^n$.

Now we assume some homogeneous structure of $h^*TS^n$ in the following way:
We assume that $K=SO(m+1)$ is acting on both $S^m$ and all 
$T_xS^m\subset\R^{m+1}$ in the standard way (and we do not distinguish between
$k$ and $k_*$ here). On $S^n$ we assume an operation of $SO(m+1)$ by some
representation $\l:SO(m+1)\to SO(n+1)$
(and hence on vectors in $E$ by the same matrices). Both operations
of $SO(m+1)$ are isometric. Moreover, we assume that
$h$ is $K$-equivariant, which in our case means
\[
  h(kx)=\l_kh(x)
\]
for all $k\in SO(m+1)$, $x\in S^m$. 

Because of
\eas
  \d_{ku}(Y(k^{-1}x))&=&(\d_uY)(k^{-1}x),\\
  \d_{ku}(\tau_kY)(x)&=&\l_k(\d_uY)(k^{-1}x),\\
  \D_{ku}(\tau_kY)(x)&=&\l_k(\D_uY)(k^{-1}x)\\
  &=&\tau_k\D_uY,
\eeas
$\D$ is $SO(m+1)$-equivariant, which means that $\D$ can be used as
reference connection as in the last section. We want to investigate for
which $A$ the connection $D=D_A=\D+A$ is Yang-Mills.
We fix $x\in S^m$, $v\in T_xS^m$, and consider 
a path in $SO(m+1)$ given by
\[
  k(t):=\id+x\otimes\{(\cos t-1)x+(\sin t)v\}
    +v\otimes\{(\cos t-1)v-(\sin t)x\}.
\]
We observe
\eas
  k(0)&=&\id,\\
  k'(0)&=&x\otimes v-v\otimes x,\\
  \l_k'(0)&=&h(x)\otimes\d_vh(x)-\d_vh(x)\otimes h(x).
\eeas
Differentiating the equivariance relation for $A$, we find
\eas
  0&=&\frac{d}{dt}_{|t=0}\Big(\l_{k(t)}^{-1}A_{k(t)u}(k(t)x)\l_{k(t)}\Big)\\
  &=&[A_u(x),\l_k'(0)]+A_{k'(0)u}(x)+A'_u(x)k'(0)x\\
  &=&[A_u(x),h(x)\otimes\d_vh(x)-\d_vh(x)\otimes h(x)]
    +A_{\<x,u\>v-\<v,u\>x}(x)+\d_vA_u(x)\\
  &=&h(x)\otimes(A_u(x)\d_vh(x))-(A_u(x)\d_vh(x))\otimes h(x)+\d_vA_u(x),
\eeas
which means that all derivatives of $A$ can be expressed by terms of order
zero:
\[
  \d_vA_u=(A_u\d_vh)\otimes h-h\otimes(A_u\d_vh).
\]
This implies
\[
  (\d_vA_u)Y=\<A_u\d_vh,Y\>h
\]
for sections $Y\in\O^0(h^*TS^n)$. We use this to calculate further
(with the first ``$=$'' being the definition of $\D$ extended to forms)
\eas
  (\D_vA_u)Y&=&\D_v(A_uY)-A_u\D_vY\\
  &=&\d_v(A_uY)+\<\d_vh,A_uY\>h-A_u\d_vY-A_u(\<\d_vh,Y\>h)\\
  &=&(\d_vA_u)Y+\<\d_vh,A_uY\>h\\
  &=&\<A_u\d_vh,Y\>h+\<\d_vh,A_uY\>h\\
  &=&0
\eeas
because $A_u\in\so(n+1)$ is skew-symmetric.
Therefore $A_u$ is covariant constant with respect to $\D$,
\[
  \D A=0.
\]
Knowing the curvature of $\D$,
\[
  (F_0)_{uv}=\d_vh\otimes\d_uh-\d_uh\otimes\d_vh,
\]
we infer that
\eas
  (F_A)_{uv}&=&(F_0)_{uv}+\D_uA_v-\D_vA_u+[A_u,A_v]\\
  &=&\d_vh\otimes\d_uh-\d_uh\otimes\d_vh+[A_u,A_v].
\eeas
Using the fact that $|F_A|^2$ is constant due to the transitivity of the
$SO(m+1)$-action on $S^m$, we conclude that (without integration) 
\[
  c\,\YM(A)=\frac12\sum_{u,v=1}^m\Big|\d_vh(e)\otimes\d_uh(e)
    -\d_uh(e)\otimes\d_vh(e)+[A_u(e),A_v(e)]\Big|^2
\]
for some dimension-dependent constant $c>0$, where we abbreviate $e=e_{m+1}$
and $\d_u=\d_{e_u}$. The first variation of this functional is (where from
now on we omit the argument $(e)$)
\[
  c\,\delta\YM(A,\Phi)
  =\sum_{u,v=1}^m\Big\<\d_vh\otimes\d_uh
    -\d_uh\otimes\d_vh+[A_u,A_v]\,,\,
    [\Phi_u,A_v]+[A_u,\Phi_v]\Big\>,
\]
from which we read off its ``Euler-Lagrange'' equation, which in this case
is just some system of algebraic equations:
\eas
  \lefteqn{\sum_{u,v=1}^m\sum_{i,j=1}^n\Big[
    \Big(\d_vh^j\d_uh^i-\d_uh^j\d_vh^i+\sum_{k=1}^n(A_u^{ik}A_v^{kj}
      -A_v^{ik}A_u^{kj})\Big)}\\
  &&\cdot\sum_{k=1}^n\Big(\Phi_u^{ik}A_v^{kj}-A_v^{ik}\Phi_u^{kj}
    +A_u^{ik}\Phi_v^{kj}-\Phi_v^{ik}A_u^{kj}\Big)\Big]\,=\,0
\eeas
for every choice of real numbers $\Phi^{ij}_u$ for $u\in\{1,\ldots,m\}$ and
$i,j\in\{1,\ldots,n\}$ satisfying $\Phi_u^{ij}=-\Phi_u^{ji}$, where we have 
assumed w.l.o.g.\ that $h(e_{m+1})=e_{n+1}$. Choosing 
$\Phi_z^{pq}=-\Phi_z^{qp}=1$ for fixed $p<q$ and $z$, and $\Phi_u^{ij}=0$ 
in all other cases, we can write this a s a system of $\frac12mn(n-1)$ cubic
equations for the same number of variables $A_u^{ij}$ (cf.\ \cite{PU} for
another formulation of these cubic equations). Therefore, in 
principle, we know ``all'' Yang-Mills connections with the symmetries 
considered. But we will not have to go into any further detail, because
all we want to know for the purpose of this paper is that $A\equiv0$ is
always a solution, which we easily read off from the equation above.
To be more precise, we have proven that $\D=D_0$ is critical for $\YM$
with respect to equivariant variations (and maybe a few more). 
But by the ``symmetric criticality''
Proposition \ref{psc}, this means that $\D$ is actually
Yang-Mills. Hence we have proven:
\begin{quote}\em
  In the special homogeneous setting considered here, the pull-back
  $\D$ of the Levi-Civita connection of $TS^n$ via $h$ is a
  Yang-Mills connection.
\end{quote}
This should have been well-known, and probably follows from the results in 
\cite{PU} or even from Itoh's earlier paper \cite{It}, and we have given the
proof mainly to introduce our setting and notation (which differ significantly
from theirs).

\strut\\
{\bf Remark.} It seems not to be true that the pull-back of the
Levi-Civita connection of $S^n$ via a homogeneous mapping $h:S^m\to S^n$
is always Yang-Mills. For example, the Hopf map $h:S^3\to S^2$ is 
$U(2)$-equivariant with $U(2)$ acting transitively on $S^3$. However,
the $\D$ that we obtain from $h$ is not Yang-Mills, as a direct
calculation shows. The same probably applies for the Hopf maps
$S^7\to S^4$ and $S^{15}\to S^8$. The only point in the above proof that
does not carry over is the choice of the path $k(t)$ in the symmetry group
$K$. For the Hopf examples, this group is no longer $SO(m+1)$, and this is
where the argument fails.

\strut\\
{\bf Examples.} Nevertheless, there are enough examples of representations
of $SO(m+1)$ to make the above considerations interesting for us:

{\bf (o)} If $m=n$ and $SO(m+1)$ acts also on the target sphere in the
standard way, then $h:S^m\to S^m$ is the identity and $\D$ is the
Levi-Civita connection of $TS^m$.

{\bf (i)} If $m=n=1$, we consider $e^{i\vartheta}\in SO(2)$ acting on
the target circle as $e^{i\ell\vartheta}$ for some $\ell\in\N$. The connection
$\D$ is the flat one, but $h:S^1\to S^1$ here is the mapping
$z\mapsto z^\ell$.

{\bf (ii)} Let $\ell\in\N$ and identify $\R^{n+1}$ with the space of 
$\ell$-homogeneous
harmonic polynomials on $R^{m+1}$, which implies $n=\frac{(2\ell+m-1)
(\ell+m-2)!}{\ell!(m-1)!}-1$. This identification is made
via an orthonormal basis $\{b_i\}_{1\le i\le n+1}$ of this
space (with respect to the scalar product $\<f,g\>=\int_{S^m}fg$). 
A representation $\l:SO(m+1)\to SO(n+1)$ is given by $\l_k f(x):=f(k^{-1}x)$,
and the corresponding $h_{m,\ell}:S^m\to S^n$ is given by
$h_{m,\ell}(x)=(b_1(x),\ldots,b_{n+1}(x))$. This map has been considered by
doCarmo and Wallach \cite{dCW}; it is a harmonic mapping as well as
(after rescaling the domain sphere suitably) a minimal immersion.\\
\strut
 
The geometries described in these examples have some more properties that will
be important if we want to take them as ``building blocks'' for the join
construction we will perform in this paper. We summarize what we need in
the following definition. The term ``Yang-Mills eigenmap'' here is
motivated to some extent by property (iii), but more by the fact that
joining {\em harmonic maps\/} is based on a similar concept called
``harmonic eigenmap''. We write $d_{LC}$ for Levi-Civita connections.

\begin{definition}[Yang-Mills eigenmap]
We call a map $h:S^m\to S^n$ a {\em Yang-Mills eigenmap\/}, if there exist
numbers $\l>0$, $\mu\ge0$, such that
\eas
  (i)&&h^*d_{LC}\mbox{ is a Yang-Mills connection,}\\
  (ii)&&|dh|^2\equiv\l,\\
  (iii)&&\sum_vF(h^*d_{LC})_{uv}\d_v h=\mu\d_u h
    \qquad\forall u\in\{1,\ldots,m\}.
\eeas
\end{definition}

Here, by $\sum_v$ and $\sum_V$, we mean summation over orthonormal
bases of $T_xS^m$ or $T_{h(x)}S^n$, respectively. 

An immediate consequence of $(iii)$ is that $|F(h^*d_{LC})|^2$ is
constant. Since 
  $$F(h^*d_{LC})_{uv}Y=\<\d_vh,Y\>\d_uh-\<\d_uh,Y\>\d_vh,$$
we have
\begin{eqnarray}\label{Fconst}
  |F(h^*d_{LC})|^2&=&\sum_{u,v}\sum_{U,V}\Big(\<\d_vh,Y\>\<\d_uh,Z\>
    -\<\d_uh,Y\>\<\d_vh,Z\>\Big)^2\nn\\
  &=&2\sum_{u,v}\Big(|\d_uh|^2|\d_vh|^2-\<\d_uh,\d_vh\>^2\Big)\nn\\
  &=&2\sum_{u,v}\<F(h^*d_{LC})_{uv}\d_vh,\d_uh\>\nn\\
  &\equiv&2\l\mu.
\end{eqnarray}

\strut\\
{\bf Examples of Yang-Mills eigenmaps\/} arise from the examples of
homogeneous Yang-Mills connections above. We have

{\bf (o)} the identities $\id_m:S^m\to S^m$, with $\l=m$, $\mu=m-1$;

{\bf (i)} the mappings $d_\ell:S^1\to S^1$, $d_\ell(z)=z^\ell$ with 
$\l=\ell^2$, $\mu=0$;

{\bf (ii)} the standard immersions 
$h_{m,\ell}:S^m\to S^{\frac{(2\ell+m-1)(\ell+m-2)!}{\ell!(m-1)!}-1}$
described above, with $\l=\ell(\ell+m-1)$ and
$\mu=\frac{m-1}m\,\ell(\ell+m-1)$.

Strictly speaking, (o) and (i) are special cases of (ii); we list them
separately because of their distinctive geometric features. 

All these examples happen to be
harmonic eigenmaps, too. Detailed accounts of harmonic eigenmaps and
related concepts can be found in the books \cite{Ba} and \cite{ER1}.
As we do not have any more examples than the ones listed here, we do not know 
whether every Yang-Mills eigenmap is automatically also a harmonic eigenmap.
Nor do we know whether we should prepare for Yang-Mills eigenmaps which
do not come from group representations (they do exist in the case of harmonic
eigenmaps). 

\begin{remark}
In all of our examples of Yang-Mills eigenmaps, we observe 
$\mu=\frac{m-1}m\,\l$. We have not assumed that in the definition of 
eigenmaps, because we will not need it in our discussion of reduction
of Yang-Mills to an o.d.e., nor is it needed for the sufficient conditions
for solving that o.d.e. The only point where it might prove important is
the question whether the conditions obtained are also necessary, cf.\
Remark \ref{sharp} below.
\end{remark}

\section{Topological motivation of our ansatz}

Now we have to justify the special kind of ansatz we are going to make
below. To this end, we write $S^{m_1+m_2+1}$ as the join $S^{m_1}*S^{m_2}$,
that is the warped product
  $$S^{m_1+m_2+1}\cong[0,\pi/2]\times_{\cos^2}S^{m_1}\times_{\sin^2}S^{m_2}$$
which closes smoothly across the endpoints of $[0,\pi/2]$. Assume we are given
an $SO(n_1)$-vectorbundle $E_1\to S^{m_1}$ and an $SO(n_2)$-vectorbundle
$E_2\to S^{m_2}$. We want to construct a {\em join\/} $E_1*E_2$ of $E_1$ and
$E_2$ as an $SO(n_1{+}n_2{+}1)$-vectorbundle over $S^{m_1}*S^{m_2}$ by roughly
``connecting every point in a fiber of $E_1$ with every point in a fiber of
$E_2$''. To make this precise, we parametrize $S^{m_1+m_2+1}$ by three
patches:
\begin{eqnarray*}
  &&\p_1:(0,\pi/2)\times S^{m_1}\times S^{m_2}\to S^{m_1+m_2+1},\\
  &&\qquad\p_1(t,x_1,x_2):=(x_1\cos t,x_2\sin t);\\
  &&\p_2:B^{m_1+1}\times S^{m_2}\to S^{m_1+m_2+1},\\
  &&\qquad\p_2(y_1,x_2):=\exp_{(0,x_2)}(y_1,0)
    =(\tf{y_1}{|y_1|}\sin|y_1|, x_2\cos|y_1|);\\
  &&\p_3:S^{m_1}\times B^{m_2+1}\to S^{m_1+m_2+1},\\
  &&\qquad\p_3(x_1,y_2):=\exp_{(x_1,0)}(0,y_2)
    =(x_1\cos|y_2|,\tf{y_2}{|y_2|}\sin|y_2|).
\end{eqnarray*}
On the image of $\p_1$ (which is all of $S^{m_1+m_2+1}$ except for two
``singular spheres''), the join of $E_1$ and $E_2$ is easily described: It is
simply the $\p_1^{-1}$-pullback of the product bundles $T\times E_1\times
E_2$, where $T$ is the trivial $\R$-bundle over $(0,\pi/2)$. The question now
is, under which condition this bundle closes smoothly across the singular
spheres to give a smooth $SO(n_1{+}n_2{+}1)$-vectorbundle over
$S^{m_1+m_2+1}$. 

This is a topological condition on the bundles $E_1$ and $E_2$, and it can be
formalized as follows. Denote by $T$ now the trivial $\R$-bundle over
$(0,1]$. Then the product bundle $T\times E_i\to (0,1]\times S^{m_i}$ must be
the pull-back of an $SO(m_i{+}1)$-bundle $\widetilde{E}_i\to B^{m_i+1}$ via the
map $(t,x_i)\mapsto tx_i$. Since every bundle over $B^{m_i+1}$ is trivial, it
is no restriction to assume
$\widetilde{E}_i=B^{m_i+1}\times\R^{n_i+1}$. Moreover, for every
$y_i\in\R^{m_i+1}$, the fiber $\R^{n_i+1}$ over $y_i$ contains a well-defined
direction which corresponds to the positive $T$-direction in the product
bundle. This defines a mapping $h_i:S^{m_i}\to S^{n_i}$ for which
$E=h_i^*TS^{n_i}$. Therefore the only bundles $E_i\to S^{m_i}$ ($i\in\{1,2\}$)
for which a smooth join can be defined are the pull-back bundles
$h_i^*TS^{n_i}$ for a pair of maps $h_i:S^{m_i}\to S^{n_i}$.

Given such a pair, we still have to find out, for which connections on
$h_i^*TS^{n_i}$ a suitable ansatz will reduce the Yang-Mills equation to an
o.d.e.\ system. 
Of course, we must think of such connections as being Yang-Mills
and ``totally homogeneous'' in a suitable sense. It turns out that suitable
``building blocks'' for our construction will be the pull-backs of the
Levi-Civita-connections of $TS^{n_i}$ under the Yang-Mills eigenmaps
defined and discussed above.

\section{Reduction}

We consider the sphere $S^{m_1+m_2+1}$ represented (somewhat sloppy concerning
the interval endpoints) as the doubly warped product
  $$(M,\gamma)
    :=[0,\pi/2]\times_{\cos^2}S^{m_1}\times_{\sin^2}S^{m_2}.$$
The Riemannian manifold $(M,\gamma)$ is isometric to the sphere with the
standard Euclidean metric.

As indicated above,
we consider an $SO(n_1+n_2+1)$-bundle over $M$ which is given as follows:
Let $h_1:S^{m_1}\to S^{n_1}$ and $h_2:S^{m_2}\to S^{n_2}$ be Yang-Mills
eigenmaps. The bundle $\Phi^*TN\to M$ under consideration is the 
pull-back of
the tangent bundle of the warped product
  $$[0,\pi/2]\times_{\cos^2}S^{n_1}\times_{\sin^2}S^{n_2}=:(N,g)
      \cong S^{n_1+n_2+1}$$
via the map
  $$\Phi:=(\id,h_1,h_2):M\to N,$$
where the parameters $(\l,\mu)$ from above
are now denoted by $(\l_i,\mu_i)$ for $h_i$.
As discussed in the previous section, this bundle can be viewed as a bundle on
all of $S^{n_1+n_2+1}$, closing smoothly across the endpoints of $[0,\pi/2]$.

The $h_i$-pullback of the Levi-Civita connection of $S^{n_i}$ will be denoted
by $\D^i$ and its curvature by $F^i$. In particular,
$F^1$ and $F^2$ are given by
\begin{eqnarray*}
  F^1_{uv}U&=&g(\d_vh_1,U)\d_uh_1-g(\d_uh_1,U)\d_vh_1,\\
  F^2_{wz}W&=&g(\d_zh_2,W)\d_wh_2-g(\d_wh_2,W)\d_zh_2.
\end{eqnarray*} 

In what follows, we denote the variable in $[0,\pi/2]$ by $t$, and the vector
field $\frac{\d}{\d t}$ by $x$ if viewed as a vector field in $TM$, and by $X$
when viewed as a vector field in $\Phi^*TN$. By $u,v$ we mean vector fields 
in $TM$ tangential to $S^{m_1}$ and by $w,z$ tangential to
$S^{m_2}$. Similarly, $U,V$ denote vector fields in $\Phi^*TN$ tangential to
$S^{n_1}$ and $W,Z$ vector fields tangential to $S^{n_2}$.

The pull-back $\D$ of the Levi-Civita connection of $N\cong S^{n_1+n_2+1}$ 
by $\Phi$ is characterized by
\begin{eqnarray*}
  \D_xX&=&0,\\
  \D_uV&=&\cos(t)^{-1}\D^1_uV+\tan(t)g(\d_uh_1,V)X,\\
  \D_wZ&=&\sin(t)^{-1}\D^2_wZ-\cot(t)g(\d_wh_2,Z)X,\\
  \D_uX&=&-\tan(t)\d_uh_1,\\
  \D_wX&=&\cot(t)\d_wh_2,\\
  \D_xV&=&-\tan(t)V,\\
  \D_xZ&=&\cot(t)Z,\\
  \D_uZ&=&0,\\
  \D_wV&=&0.
\end{eqnarray*}
The ansatz we make for our connection $D$ on $\Phi^*TN$ differs only slightly
from that, in an ``equivariant'' way:
\begin{eqnarray*}
  D_xX&=&0,\\
  D_uV&=&\frac1{\cos(t)}\,(\D^1_uV+\a(t)g(\d_uh_1,V)X),\\
  D_wZ&=&\frac1{\sin(t)}\,(\D^2_wZ-\b(t)g(\d_wh_2,Z)X),\\
  D_uX&=&-\frac{\a(t)}{\cos(t)}\,\d_uh_1,\\
  D_wX&=&\frac{\b(t)}{\sin(t)}\,\d_wh_2,\\
  D_xV&=&-\tan(t)V,\\
  D_xZ&=&\cot(t)Z,\\
  D_uZ&=&0,\\
  D_wV&=&0.
\end{eqnarray*}
This connection is still metric with respect to $g$. The basic idea for
finding $\a$ and $\b$ for which $D$ is Yang-Mills will be minimizing
the Yang-Mills functional over $SO(m_1+1)\times SO(m_2+1)$-equivariant
connections (which $D$ is).

Now we are ready to calculate the curvature of $D$, which we denote by $F$.
Since we know that $F$ is a tensor, i.e.\ a differential operator of order
$0$, we can assume we are calculating everything in a point where 
$\D^1_uU=0$ etc.\ and $\D^1_u\d_vh_1=\D^1_v\d_uh_1$ etc.:
\begin{eqnarray*}
  F_{uv}U&=&(D_uD_v-D_vD_u)U\\
  &=&D_u(\tf1{\cos}\D^1_vU+\tf{\a}{\cos}g(\d_vh_1,U)X)
      -D_v(\tf1{\cos}\D^1_uU+\tf{\a}{\cos}g(\d_uh_1,U)X)\\
  &=&\tf1{\cos^2}(\D^1_u\D^1_vU-\D^1_v\D^1_uU)\\
  &&{}+\tf{\a}{\cos}(g(\D^1_u\d_vh_1-\D^1_v\d_uh_1,U)
    +g(\d_vh_1,\D^1_uU)-g(\d_uh_1,\D^1_vU))X\\
  &&+\tf{\a^2}{\cos^2}(g(\d_uh_1,\d_vh_1)g(X,U)-g(\d_vh_1,\d_uh_1)g(X,U))X\\
  &&+\tf{\a^2}{\cos^2}(g(\d_uh_1,U)g(\d_vh_1,X)-g(\d_vh_1,U)g(\d_uh_1,X))X\\
  &&+\tf{\a^2}{\cos^2}(-g(\d_vh_1,U)\d_uh_1+g(\d_uh_1,U)\d_vh_1)\\
  &=&\tf1{\cos^2}F^1_{uv}U
    +\tf{\a^2}{\cos^2}(g(\d_uh_1,U)\d_vh_1-g(\d_vh_1,U)\d_uh_1)\\
  &=&\tf{\a^2-1}{\cos^2}(g(\d_uh_1,U)\d_vh_1-g(\d_vh_1,U)\d_uh_1),\\
  F_{uv}W&=&0,\\
  F_{uv}X&=&-D_u(\tf{\a}{\cos}\d_vh_1)+D_v(\tf{\a}{\cos}\d_uh_1)\\
  &=&-\D^1_u(\tf{\a}{\cos^2}\d_vh_1)+\D^1_v(\tf{\a}{\cos^2}\d_uh_1)
    -\tf{\a^2}{\cos^2}(g(\d_uh_1,\d_vh_1)-g(\d_vh_1,\d_uh_1))X\\
  &=&0,\\
  F_{wz}W&=&\tf{\b^2-1}{\sin^2}(g(\d_wh_1,W)\d_zh_1-g(\d_zh_1,W)\d_wh_1),\\
  F_{wz}U&=&0,\\
  F_{wz}X&=&0,\\
  F_{xu}U&=&D_x(\tf1{\cos}\D^1_uU+\tf{\a}{\cos}g(\d_uh_1,U)X)-D_u(-\tan U)\\
  &=&\tf{\a'+\a\tan}{\cos}g(\d_uh_1,U)X-2\tf{\a\tan}{\cos} 
    g(\d_uh_1,U)X+\tf{\a\tan}{\cos}g(\d_uh_1,U)X\\
  &=&\tf{\a'}{\cos}g(\d_uh_1,U)X\\
  F_{xu}W&=&D_x(D_uW)-D_u(\cot W)\\
  &=&0,\\
  F_{xu}X&=&D_xD_uX\\
  &=&-D_x(\tf{\a}{\cos}\d_uh_1)\\
  &=&-(\tf{\a'+\a\tan}{\cos}-\tf{\a\tan}{\cos})\d_uh_1\\
  &=&-\tf{\a'}{\cos}\d_uh_1,\\
  F_{xw}U&=&0,\\
  F_{xw}W&=&-\tf{\b'}{\sin}g(\d_wh_2,W)X,\\
  F_{xw}X&=&-\tf{\b'}{\sin}\d_wh_2,\\
  F_{uw}U&=&-D_wD_uU\\
  &=&-D_w(\tf1{\cos}\D^1_uU+\tf{\a}{\cos}g(\d_uh_1,U)X)\\
  &=&-\tf{\a}{\cos}g(\d_uh_1,U)D_wX\\
  &=&-\tf{\a\b}{\cos\sin}g(\d_uh_1,U)\d_wh_2,\\
  F_{uw}W&=&\tf{\a\b}{\cos\sin}g(\d_wh_2,W)\d_uh_1,\\
  F_{uw}X&=&D_u(\cot\d_wh_2)-D_w(\tan\d_uh_1)\\
  &=&0.
\end{eqnarray*}
Summing over $\gamma$-orthonormal bases $u,v,w,z$ and $g$-orthonormal bases
$U,W$, we infer
  $$|F|^2=4\l_1\,\frac{\a'^2}{\cos^2}+4\l_2\frac{\b'^2}{\sin^2}
    +2\l_1\mu_1\,\frac{(\a^2-1)^2}{\cos^4}
    +2\l_2\mu_2\,\frac{(\b^2-1)^2}{\sin^4}
    +4\l_1\l_2\frac{\a^2\b^2}{\cos^2\sin^2}.$$
Therefore, up to a constant depending only on $m_1$, $m_2$, the Yang-Mills
functional of $D$ equals
\begin{eqnarray}\label{j}
  J(\a,\b)&:=&\int_0^{\pi/2}
    \Big\{\frac{2\l_1}{\cos^2}\,\a'^2
    +\frac{2\l_2}{\sin^2}\,\b'^2+\frac{2\l_1\l_2}{\cos^2\sin^2}\,\a^2\b^2
    \nonumber\\
  &&\quad{}+\frac{\l_1\mu_1}{\cos^4}(\a^2-1)^2
    +\frac{\l_2\mu_2}{\sin^4}(\b^2-1)^2
    \Big\}\cos^{m_1}\sin^{m_2}\,dt.
\end{eqnarray}
The Euler-Lagrange equations of $J$ are
\begin{eqnarray}\label{ela}
  \a''+(m_2\cot-(m_1-2)\tan)\a'
    -\frac{\mu_1}{\cos^2}\,(\a^3-\a)
    -\frac{\l_2}{\sin^2}\,\a\b^2&=&0,\\
\label{elb}
  \b''+((m_2-2)\cot-m_1\tan)\b'
    -\frac{\mu_2}{\sin^2}\,(\b^3-\b)
    -\frac{\l_1}{\cos^2}\,\a^2\b&=&0.
\end{eqnarray}
The reduction setting is made in such a way that stationary points of the
reduced functional $J$ represent Yang-Mills fields:

\begin{proposition}[reduction theorem]
The connection $D$ is a smooth Yang-Mills connection on $\Phi^*TN$ 
if and only if the functions
$\a,\b:[0,\pi/2]\to\R$ are solutions of \re{ela}, \re{elb}
with the boundary values
\begin{equation}\label{bc}
  \a(0)=0,\qquad\a(\pi/2)=1,\qquad\b(0)=1,\qquad\b(\pi/2)=0.
\end{equation}
\end{proposition}

\proof
To calculate the Yang-Mills equations for our setting, we have to
differentiate $F$. In the following calculations, we make the same assumptions
on the vector fields as above. We also use the fact that $F^1$, $F^2$ are
Yang-Mills connections, and assume summation if an ``index'' is repeated.
\begin{eqnarray*}
  (D_uF_{uv})U&=&D_u(F_{uv}U)-F_{uv}D_uU\\
  &=&\tf{\a^2-1}{\cos^2}D_u(g(\d_uh_1,U)\d_vh_1-g(\d_vh_1,U)\d_uh_1)\\
  &&{}-\tf1{\cos}F_{uv}\D^1_uU-\tf{\a}{\cos}g(\d_uh_1,U)F_{uv}X\\
  &=&({\nabla^1}^*F^1)_vU+\tf{\a^3-\a}{\cos^3}(g(\d_uh_1,U)g(\d_uh_1,\d_vh_1)\\
  &&{}-g(\d_vh_1,U)g(\d_uh_1,\d_uh_1))X,\\
  &=&\tf{\a^3-\a}{\cos^3}(g(\d_uh_1,U)g(\d_uh_1,\d_vh_1)
    -g(\d_vh_1,U)g(\d_uh_1,\d_uh_1))X,\\
  (D_uF_{uv})W&=&0,\\
  (D_uF_{uv})X&=&D_u(F_{uv}X)+F_{uv}(\tf{\a}{\cos}\d_uh_1)\\
  &=&\tf{\a^3-\a}{\cos^3}(g(\d_uh_1,\d_uh_1)\d_vh_1
    -g(\d_uh_1,\d_vh_1)\d_uh_1),\\
  (D_wF_{wv})U&=&D_w(\tf{\a\b}{\cos\sin}g(\d_vh_1,U)\d_wh_2)\\
  &=&-\tf{\a\b^2}{\cos\sin^2}g(\d_vh_1,U)g(\d_wh_2,\d_wh_2)X,\\
  (D_wF_{wv})W&=&-D_w(\tf{\a\b}{\sin\cos}g(\d_wh_2,W)\d_vh_1)\\
  &=&0,\\
  (D_wF_{wv})X&=&-F_{wv}D_wX\\
  &=&-\tf{\b}{\sin}F_{wv}\d_wh_2\\
  &=&\tf{\a\b^2}{\cos\sin^2}g(\d_wh_2,\d_wh_2)\d_vh_1,\\
  (D_xF_{xv})U&=&D_x(\tf{\a'}{\cos}g(\d_vh_1,U)X)-F_{xv}(-\tan U)\\
  &=&(\tf{\a''+\a'\tan}{\cos}-2\tf{\a'\tan}{\cos}
    +\tf{\a'\tan}{\cos})g(\d_vh_1,U)X\\
  &=&\tf{\a''}{\cos}g(\d_vh_1,U)X,\\
  (D_xF_{xv})W&=&0,\\
  (D_xF_{xv})X&=&-D_x(\tf{\a'}{\cos}\d_vh_1)\\
  &=&(-\tf{\a''+\a'\tan}{\cos}+\tf{\a'\tan}{\cos})\d_vh_1\\
  &=&-\tf{\a''}{\cos}\d_vh_1;\\
  D_uF_{ux}&=&0;
\end{eqnarray*}
and similar terms for $S^{n_2}$-components.

The next thing we have to check is what $D^*$ looks like in our
coordinates. There are induced metrics from $g$ for $1$-forms and
$2$-forms, which we again denote by $g$. 
By the definition of $D^*$ and partial integration, we find
for every $2$-form $G$ and every $1$-form~$\p$
  \def\iiint{\int_0^{2\pi}\int_{S^{m_1}}\int_{S^{m_2}}}
  \def\dvvv{\,d\mbox{\rm vol}_2\,d\mbox{\rm vol}_1\,dt}
\begin{eqnarray*}
  -\lefteqn{\iiint g(D^*G,\p)\cos^{m_1}\sin^{m_2}\dvvv}\\
  &=&-\iiint g(G,D\p)\cos^{m_1}\sin^{m_2}\dvvv\\
  &=&\iiint g(D\cdot G,\p)\cos^{m_1}\sin^{m_2}\dvvv\\
  &&+\iiint g(G(x,\,\cdot\,),\p)(m_2\cos^{m_1+1}\sin^{m_2-1}
      -m_1\cos^{m_1-1}\sin^{m_2+1})\dvvv\\
  &&+\iiint (\partial_tg)(G(x,\,\cdot\,),\p)
      \cos^{m_1}\sin^{m_2}\dvvv.
\end{eqnarray*}
Knowing that
\[
  \d_tg=2\diag(-(\tan t)\id_{m_1},(\cot t)\id_{m_2},0)\,g,
\]
we can read off from the previous equation how $D^*$ operates.
This is combined with the calculation above and (ii), (iii) to give
\begin{eqnarray*}
  -(D^*F)_vU&=&(D\cdot F)_vU+(m_2\cot-(m_1-2)\tan)F_{xv}U\\
    &=&\Big\{-\frac{\mu_1}{\cos^3}(\a^3-\a)-\frac{\l_2}{\cos\sin^2}\a\b^2
        +\frac{\a''}{\cos}\\
    &&\quad{}+\frac{\a'}{\cos}(m_2\cot-(m_1-2)\tan)\Big\}g(\d_vh_1,U)X,\\
  -(D^*F)_vX&=&(D\cdot F)_vX+(m_2\cot-(m_1-2)\tan)F_{xv}X\\
    &=&\Big\{\frac{\mu_1}{\cos^3}(\a^3-\a)+\frac{\l_2}{\cos\sin^2}\a\b^2
        -\frac{\a''}{\cos}\\
    &&\quad{}-\frac{\a'}{\cos}(m_2\cot-(m_1-2)\tan)\Big\}\d_vh_1,
\end{eqnarray*}
and the corresponding
equations for the $S^{n_2}$ components. (Some components always vanish.)
This proves that $D$
is Yang-Mills away from $t\in\{0,\pi/2\}$ if and only if \re{ela} and \re{elb}
are satisfied.

Now we turn to the boundary conditions. Because $\a=\sin$ and $\b=\cos$
correspond to the pullback of the Levi-Civita connection of $TS^{n_1+n_2+1}$,
the boundary conditions \re{bc} make sure that the connection is continuous
even across the singular orbits $\{t=0\}$ and $\{t=\pi/2\}$. For $D$ to be
of class $C^1$, $\a$ and $\b$ also have to satisfy
\begin{equation}\label{bc2}
  \a'(0)=1,\qquad\a'(\pi/2)=0,\qquad\b'(0)=0,\qquad\b'(\pi/2)=-1.
\end{equation}
But this is easily seen to hold for any solution of the boundary value problem
made of \re{ela}, \re{elb}, \re{bc}. Once this is checked, the parity of the
differential equations \re{ela} and \re{elb} implies that for any solution
with the boundary values \re{bc} and \re{bc2} the function $\a$ is odd with
respect to $t=0$ and even with respect to $t=\pi/2$, while for $\b$ the
opposite holds. But those are exactly the conditions to ensure that $D$ is
smooth across the singular orbits. This proves the reduction theorem.\qed

From the harmonic map analogon of our problem, we know that the substitution
  $$\a(t)=A(\log(\tan t)),\qquad\b(t)=B(\log(\tan t))$$
is useful. With $s=\log(\tan t)$ we calculate
\begin{eqnarray*}
  \a'(t)&=&(e^s+e^{-s})A'(s),\\
  \a''(t)&=&(e^s+e^{-s})^2A''(s)+(e^{2s}-e^{-2s})A'(s)
\end{eqnarray*}
etc., which transforms \re{ela}--\re{bc} to give
\begin{eqnarray}\label{eta}
  A''-\frac{(m_1-3)e^s-(m_2-1)e^{-s}}{e^s+e^{-s}}\,A'
    -\frac{\mu_1e^s}{e^s+e^{-s}}\,(A^3-A)
    -\frac{\l_2 e^{-s}}{e^s+e^{-s}}\,AB^2&=&0,\\
\label{etb}
  B''-\frac{(m_1-1)e^s-(m_2-3)e^{-s}}{e^s+e^{-s}}\,B'
    -\frac{\mu_2e^{-s}}{e^s+e^{-s}}\,(B^3-B)
    -\frac{\l_1 e^s}{e^s+e^{-s}}\,A^2B&=&0.
\end{eqnarray}
with the boundary conditions
\begin{equation}\label{bct}
  A(-\infty)=0,\qquad A(\infty)=1,\qquad B(-\infty)=1,\qquad B(\infty)=0.
\end{equation}

\section{Existence of solutions}

\subsection{The case $\mu_1,\mu_2>0$}

The first case we consider is the case where none of the eigenconnections
which are joined is flat. This is the case $\mu_1,\mu_2>0$ which is only
possible if $m_1,m_2\ge2$.
Then we find a minimizer of $J$ by the direct
method of the calculus of variations.

\begin{lemma}[existence of minimizers]\label{exist}
Assume $m_1,m_2\ge2$ 
and $\mu_1,\mu_2>0$. Then there is a solution $(a,b)$ of
\re{ela}, \re{elb} on $(0,\pi/2)$ which minimizes $J$ among all $(\a,\b)\in
C^1((0,\pi/2))^2$. It satisfies $0\le a\le 1$ and $0\le b\le 1$.
\end{lemma}

\proof
First we write $D=\D+\om$ with a matrix-valued one-form 
$\om$ which is given by
\begin{eqnarray*}
  \om(u)&=&\frac{\a-\sin}{\cos}(\d_uh_1\otimes X-X\otimes\d_uh_1),\\
  \om(w)&=&\frac{\cos-\b}{\sin}(\d_wh_2\otimes X-X\otimes\d_wh_2),\\
  \om(x)&=&0.
\end{eqnarray*}
We find
\begin{eqnarray*}
  |\om|^2&=&\frac{2\l_1}{\cos^2}\,(\a-\sin)^2
     +\frac{2\l_2}{\sin^2}\,(\b-\cos)^2,\\
  |\D\om|^2&=&\frac{2\l_1^2}{\cos^4}\,\a^2(\a-\sin)^2
     +\frac{2\l_2^2}{\sin^4}\,\b^2(\b-\cos)^2
     +\frac{2\l_1}{\cos^2}\,(\a'-\cos)^2
     +\frac{2\l_2}{\sin^2}\,(\b'+\sin)^2.
\end{eqnarray*}
It is easily read off from this that
\begin{equation}\label{coer}
  \int_M(|\om|^4+|\D\om|^2)\le c(1+J(\a,\b))
\end{equation}
with a constant $c$ depending on $\l_1,\l_2,\mu_1,\mu_2$.
Now we consider a minimizing sequence of connections $D$ with $\om\in
W^{1,2}\cap L^4$ of the special form 
considered in this paper for the Yang-Mills functional $\YM$. Since $\YM$ is
lower semi-continuous on $W^{1,2}\cap L^4$, and since we have just checked that
the $W^{1,2}$ and $L^4$ norms of $\om$ stay bounded for such a sequence, there
is a connection minimizing $\YM$ of the form considered in $W^{1,2}\cap L^4$,
which must be continuous in the orbits over $(0,\pi/2)$. This connection
is represented by a minimizer $(a,b)$ of $J$. Since minimizers of $J$ satisfy
its Euler-Lagrange equations, $(a,b)$ must be smooth on $(0,\pi/2)$.

The next step is to prove that there are minimizers with values in $[0,1]$.
To this end, consider
$f:[0,\pi/2]\to\R$ and define
  $$\tilde{f}(x):=\left\{\begin{array}{ll}
      |f(x)| & \mbox{ \ if }|f(x)|\le 1,\\
      \frac{1}{|f(x)|} & \mbox{ \ if }|f(x)|>1.
    \end{array}\right.$$
Then $\tilde{f}'^2\le f'^2$ and $(\tilde{f}^2-1)^2\le (f^2-1)^2$. 
Hence, if $(a,b)$ is minimizing, so is $(\tilde{a},\tilde{b})$, which means we
have found the minimizing solution stated in the lemma.\qed

\begin{remark}\label{const}
Assume $m_1,m_2\ge2$ and $\mu_1,\mu_2>0$.
There are exactly three constant solutions of \re{ela}, \re{elb} with values in
$[0,1]$, namely $(\a,\b)\equiv (0,0)$, $\equiv (0,1)$ or $\equiv (1,0)$. The
constant solution $(0,0)$ is never minimizing because of $J(0,0)>J(0,1)$ and
$J(0,0)>J(1,0)$.
\end{remark}

\begin{lemma}[nonconstant minimizers I]\label{cnc}
Assume $m_1,m_2\ge2$ and $\mu_1,\mu_2>0$.
No nonconstant $J$-minimizing solution $(\a,\b)$
of \re{ela}, \re{elb} with values in $[0,1]$ assumes the values
$\a(t)\in\{0,1\}$ or $\b(t)\in\{0,1\}$ at any $t\in(0,\pi/2)$.
\end{lemma}

\proof 
Assume $\a(t)=1$, then we have $\a'(t)=0$, which by way of \re{ela} implies
$\a''(t)=\l_2\sin^{-2}(t)\b(t)^2\ge0$. Since $\a$ has a maximum at $t$, this
can be true only if $\a''(t)=0$ and $\b(t)=0$. The latter implies $\b'(t)=0$.
By uniqueness of the solution of the boundary value problem with
$\a,\a',\b,\b'$ prescribed at $t$, we would have $(\a,\b)\equiv(1,0)$.

By the same reasoning, $\b(t)=1$ implies $(\a,\b)\equiv(0,1)$.

Now assume $\a(t)=0$, then $\a'(t)=0$, and \re{ela} implies $\a''(t)=0$.
Differentiate \re{ela} and find $\a^{(k)}(t)=0$ for all $k\in\N$. By 
the analyticity
of $\a$, we find $\a\equiv 0$. Once we have this, we find
  $$J(\a,\b)=\int_0^{\pi/2}\Big\{\frac{2\l_2}{\sin^2}\,\b'^2+
    \frac{\l_2\mu_2}{\sin^4}\,(\b^2-1)^2+\frac{\l_1\mu_1}{\cos^4}
    \Big\}\cos^{m_1}\sin^{m_2}\,dt,$$
which is infinity in case $m_1\le3$ (a contradiction) or easily seen to be
minimized by $\b\equiv 1$.

By the same reasoning, $\b(t)=0$ gives a contradiction or $(\a,\b)\equiv
(1,0)$.\qed

\begin{lemma}[nonconstant minimizers II]\label{nonconst}
Assume $m_1,m_2\ge2$ and $\mu_1,\mu_2>0$. If neither of the constant 
solutions $(0,1)$ and $(1,0)$ is minimizing, there is a solution of the
boundary value problem \re{ela}--\re{bc}.
\end{lemma}

\proof
Since $(0,0)$ is never minimizing, the assumption implies that
there has to be a nonconstant minimizer $(\a,\b)$ of $J$ with values in
$[0,1]$. 

In case $m_1\in\{2,3\}$ we see that $J(\a,\b)<\infty$ only if
$\a(\pi/2)=1$. But then \re{elb} implies that $\b(\pi/2)=0$ (using the fact
that $\b'$
cannot explode like $\frac1{t-\pi/2}$ if $J$ is finite).

Similarly, $m_2\in\{2,3\}$ implies $\b(0)=1$ and $\a(0)=0$.

No we consider $m_2\ge 4$. Observe that the only boundary values for $\b(0)$
that the equations \re{ela}, \re{elb} allow are $-1$, $0$ or $1$. We have
already ruled out $-1$. If $\b(0)$ was $0$, the corresponding solution $B$ of
\re{eta}, \re{etb} would asymptotically (as $s\to-\infty$) 
satisfy the linearized version of \re{eta},
  $B''+(m_2-3)B'+\mu_2B=0.$
But a fundamental system for this linearized equation consists of
$\exp((-\frac{m_2-3}2\pm\frac12\sqrt{(m_2-3)^2-4\mu_2})s)$ neither of which
is bounded at $-\infty$. Hence $B(-\infty)=0$ is not possible. This proves
$\b(0)=1$, and $\a(0)=0$ follows as in the case $m_2\le3$.

The same way we see that $\a(\pi/2)=1$ and $\b(\pi/2)=0$ also in the case
$m_1\ge4$.\qed

The lemma shows that it helps to know if the constant solutions are
minimizing. Now they are clearly not minimizing if they are unstable (in the
sense of negative directions for the second variation) or have infinite 
$J$-energy.

\begin{lemma}[unstable constant solution]\label{unst}
Assume $m_1,m_2\ge2$, $\mu_1,\mu_2>0$.
The constant solution $(\a,\b)\equiv(0,1)$ is unstable or has infinite 
$J$-energy iff
\begin{eqnarray*}
  &&m_1\in\{2,3\}\\
  \mbox{or}&&(m_1-3)^2<4\mu_1\\
  \mbox{or}&&\sqrt{(m_2-1)^2+4\l_2}\,+\sqrt{(m_1-3)^2-4\mu_1}\,<m_1+m_2-4.
\end{eqnarray*}
\end{lemma}

\proof
The case $m_1\in\{2,3\}$ is trivial. For $m_1\ge4$, we calculate
\begin{eqnarray*}
  \frac{d^2}{ds^2}_{|s=0}J(s\p,1+s\psi)
  &=&\int_0^{\pi/2}\Big\{\frac{4\l_1}{\cos^2}\,\p'^2
    +\frac{4\l_2}{\sin^2}\,\psi'^2
    +\frac{4\l_1\l_2}{\cos^2\sin^2}\,\p^2\\
  &&\qquad{}-\frac{4\l_1\mu_1}{\cos^4}\,\p^2+\frac{8\l_2\mu_2}{\sin^4}\,\psi^2
    \Big\}\cos^{m_1}\sin^{m_2}\,dt.
\end{eqnarray*}
This means that $(0,1)$ is unstable iff the quadratic form
  $$H(\p):=\int_0^{\pi/2}\Big\{\p'^2
    +\Big(\frac{\l_2}{\sin^2}-\frac{\mu_1}{\cos^2}\Big)\p^2\Big\}
    \cos^{m_1-2}\sin^{m_2}\,dt$$
becomes negative for some (bounded) function $\p$. This has been discussed by
Ding in \cite{Di} for the same function $H$ that arises with different 
constants in the construction of harmonic maps as joins of 
harmonic eigenmaps. A detailed discussion can be found in 
\cite[IX (4.4)--(4.16)]{ER2}. It shows that $H(\p)$ attains negative values if
and only if one of the three assumptions of the lemma is fulfilled.\qed  

Combining the four Lemmas from this section and the analogon of Lemma
\re{unst} for the constant solution $(1,0)$, we get our main theorem:

\begin{theorem}\label{main}
Assume $m_1,m_2\ge2$ and $\mu_1,\mu_2>0$. There is a
Yang-Mills connection of $\Phi^*TN$ corresponding to a solution
$(\a,\b)$ of the boundary value problem \re{ela}--\re{bc} if the following
conditions hold:
\begin{eqnarray*}
 (D1)\qquad&&m_1\in\{2,3\}\\
  \mbox{or}&&(m_1-3)^2<4\mu_1\\
  \mbox{or}&&\sqrt{(m_2-1)^2+4\l_2}\,+\sqrt{(m_1-3)^2-4\mu_1}\,<m_1+m_2-4
\end{eqnarray*}
and
\begin{eqnarray*}
 (D2)\qquad&&m_2\in\{2,3\}\\
  \mbox{or}&&(m_2-3)^2<4\mu_2\\
  \mbox{or}&&\sqrt{(m_1-1)^2+4\l_1}\,+\sqrt{(m_2-3)^2-4\mu_2}\,<m_1+m_2-4.
\end{eqnarray*}
\end{theorem}

\begin{remark}\label{sharp}
If $m_1=m_2$, $\l_1=\l_2$, and $\mu_1=\mu_2$, the boundary value problem
\re{ela}--\re{bc} has always a solution. This can be seen by modifying the
proof in such a way that one only minimizes over $(\a,\b)$ satisfying
$\a(\frac{\pi}2-t)=\b(t)$ for all $t\in(0,\frac{\pi}2)$. Therefore, the 
conditions $(D1)$ and $(D2)$ can only be sharp if they are automatically
satisfied in the case $m_1=m_2$, $\l_1=\l_2$, $\mu_1=\mu_2$. Unfortunately,
this is not the case if $\mu_1$ is small compared to $\l_1$. Therefore, in
the general setting considered here, the condition of Theorem \ref{main}
is sufficient, but not necessary. 

However, in our examples we always have $\mu_i=\frac{m_i-1}{m_i}\,\l_i$;
and maybe this is so for all Yang-Mills eigenmaps. Under this additional
assumption (which is the only case relevant for our construction at the
moment), $(D1)$ and $(D2)$ are always fulfilled if $m_1=m_2$, $\l_1=\l_2$,
and it is still possible that Theorem \ref{main} is sharp. We tend to expect
this to hold, because of the close analogy to the harmonic map case, where
very similar conditions have been proven to be sharp \cite{Di} \cite{PR}.
\end{remark}

\subsection{The case $m_2=1$}
An interesting case is $m_2=1$, $h_2:S^1\to S^1$ with $h_2(z)=z^k$ for some
$k\in\Z\setminus\{0\}$. 
The function $h_2$ is a Yang-Mills eigenmap, but we have $\mu_2=0$,
which means the techniques of the previous section do not apply. However, in
this particular case, we can modify the proof of 
Theorem \ref{main} to make it still
work (which we cannot do if $\D^2$ is a flat connection on a higher-dimensional
sphere). The existence theorem here reads as follows (with $\l_2=k^2$)

\begin{theorem}\label{main2} 
Assume $m_1\ge2$, $\mu_1>0$, $m_2=1$, $\mu_2=0$, $\l_2=k^2\ne0$. There is a
Yang-Mills connection of $\Phi^*TN$ corresponding to a solution
$(\a,\b)$ of the boundary value problem \re{ela}--\re{bc} if (D1) holds, which
now reads
\begin{eqnarray*}
  &&m_1\in\{2,3\}\\
  \mbox{or}&&(m_1-3)^2<4\mu_1\\
  \mbox{or}&&2\,|k|+\sqrt{(m_1-3)^2-4\mu_1}\,<m_1-3.
\end{eqnarray*}
\end{theorem}

\proof
Few changes are to be made compared to the proof of Theorem \ref{main}.
The main problem is that $J(\a,\b)$ does no longer contain a
$(\b^2-1)^2$-term. This means that $\b(0)\in\{-1,0,1\}$ is no longer needed 
to make $J$ finite. Now $\b(0)$ can take any value in $\R$ (probably), and it
will not be true that a minimizer of $J$ will more or less automatically
satisfy $\b(0)=1$. But here we can impose the boundary value $\b(0)=1$ and
minimize under this condition. 

To prove this assertion, let $(\a_n,\b_n)_{n\in\N}$ be a minimizing sequence
for $J$ under the additional hypothesis $\b(0)=1$. Again we may assume that
the images of $\a_n$ and $\b_n$ are contained in $[0,1]$. As above, we assume
that the minimizing sequence belongs to a form $\om\in W^{1,2}\cap
L^4$. Finiteness of the norms implies that $\b_n'(0)=-1$ for all $n\in\N$.

In the proof of Lemma \re{exist}, it is not immediately clear that \re{coer}
still holds, because this time the $(\b-\cos)^4$-term and the 
$\b^2(\b-\cos)^2$-Term cannot be estimated
easily by the $(\b^2-1)^2$-term which is no longer in $J$. But in the case
$m_2=1$, it can be estimated by the $\b'^2$-term of $J$ instead. This can be
seen as follows. We assume that the image of $\b$ is contained in $[0,1]$ and
that $\b(0)=1$, $\b'(0)=0$. We combine
\begin{eqnarray*}
  \int_0^1(\b(t)-\cos(t))^2\cos(t)^{m_1}\sin(t)^{-3}\,dt
  &\le&c+c\int_0^1(\b(t)-1)^2t^{-3}\,dt\\
  &\le&c+c\int_0^1\Big(\int_0^t\b'(\tau)\,d\tau\Big)^2t^{-3}\,dt\\
  &\le&c+c\int_0^1\int_0^t\b'(\tau)^2\,d\tau\,t^{-2}\,dt\\
  &\le&c+c\int_0^1\b'(t)^2\int_t^1\tau^{-2}\,d\tau\,dt\\
  &=&c+\frac{c}{2}\int_0^1\b'(t)^2(t^{-1}-1)\,dt\\
  &\le&c+c\int_0^1\b'(t)^2\cos(t)^{m_1}\sin(t)^{-1}\,dt
\end{eqnarray*}
with
\begin{eqnarray*}
  \int_1^{\pi/2}(\b-\cos)^2\cos^{m_1}\sin^{-3}\,dt
  &\le&c\int_1^{\pi/2}(\b^2+\cos^2)\cos^{m_1}\sin^{-3}\,dt\\
  &\le&c+c\int_1^{\pi/2}\b^2\cos^{m_1}\sin^{-3}\,dt\\
  &\le&c+c\int_1^{\pi/2}\b'^2\cos^{m_1}\sin^{-1}\,dt
\end{eqnarray*}
to find
\begin{equation}\label{sob}
  \int_0^{\pi/2}(\b-\cos)^2\cos^{m_1}\sin^{-3}\,dt
  \le c+c\int_0^{\pi/2}\b'^2\cos^{m_1}\sin^{-1}\,dt.
\end{equation}
Here, we have used $\b\in[0,1]$ several times, and the Sobolev inequality for
the manifold $([0,\pi/2],\cos^{m_1}\sin^{-3})$ in the second estimate. All
integrals are finite because of $\b(0)=1$ and $\b'(0)=0$. Using again
$\b\in[0,1]$, we see that \re{sob} implies
\[
  \int_0^{\pi/2}\{\b^2(\b-\cos)^2+(\b-\cos)^4\}\cos^{m_1}\sin^{-3}\,dt
  \le\int_0^{\pi/2}\b'^2\cos^{m_1}\sin^{-1}\,dt,
\]
which is exactly what was missing in the proof that \re{coer} still holds.

Once we have \re{coer}, we know that our minimizing sequence
$(\a_n,\b_n)_{n\in\N}$ stays bounded in $W^{1,2}\cap L^4$ and hence has a
weakly convergent subsequence. Again we use lower semi-continuity of $\YM$ to
conclude convergence of $(\a_n,\b_n)$ to a minimizer $(\a,\b)$ of $J$ under the
additional condition $\b(0)=1$. The only thing that remains to be checked is
that the latter condition is actually preserved in the limit. Assume it
is not. Then there is $\eps>0$ and a subsequence of $(\b_n)$, again denoted by
$(\b_n)$, such that $\min_{[0,1/n]}\b_n(x)\le 1-\eps$. Since also $\b_n(0)=1$,
this would imply
  $\int_0^{1/n}\b_n'^2\,dt\ge\frac1n\,(n\eps)^2=n\eps^2,$
which would mean $J(\a_n,\b_n)\to\infty$, a contradiction.

We have now proved, that there exist a minimizer $(\a,\b)$ of $J$ under the
additional condition $\b(0)=1$. From here, we proceed as in the proof of
Theorem \ref{main} to prove Theorem \ref{main2}.\qed

\subsection{The case $m_2=0$: Suspensions}

The case $m_2=0$ makes sense, not only formally. Remembering that
$S^0=\{-1,1\}$, we see that the join of $S^{m_1}$ and $S^0$ is 
nothing else than
the suspension $S^{m_1+1}$ of $S^{m_1}$. Consequently, we speak of Yang-Mills
suspensions here rather than of joins. This corresponds to a (simply) warped
product $M=[-\pi/2,\pi/2]\times_{\cos^2} S^{m_1}\cong S^{m_1+1}$ where we are
looking for connections on some bundle $\Phi^*TN$ of the form
\begin{eqnarray*}
  D_xX&=&0,\\
  D_uV&=&\frac1{\cos(t)}\,(\D^1_uV+\a(t)g(\d_uh_1,V)X),\\
  D_uX&=&-\frac{\a(t)}{\cos(t)}\,\d_uh_1,\\
  D_xV&=&-\tan(t)V.
\end{eqnarray*}
All notation that is not declared has a similar meaning as before.

The reduced Yang-Mills functional is
  $$J(\a):=\int_{-\pi/2}^{\pi/2}\Big\{\frac{2\l_1}{\cos^2}\,\a'^2
      +\frac{\l_1\mu_1}{\cos^4}(\a^2-1)^2\Big\}\cos^{m_1}\,dx$$
and the system of Euler-Lagrange equations reduces to a single equation
\begin{equation}\label{sus}
  \a''(t)-(m_1-2)\tan(t)\a'(t)
    -\frac{\mu_1}{\cos(t)^2}\,(\a(t)^3-\a(t))=0,
\end{equation}
with the natural boundary conditions
\begin{equation}\label{sbc}
  \a(-\pi/2)=-1,\qquad \a(\pi/2)=1.
\end{equation}
A very similar boundary value problem has been solved for harmonic suspensions
by Eells and Ratto \cite{ER1}. Since here we are very close to the harmonic
map case, we can omit details of the proof.

\begin{theorem}\label{mainsus}
Assume $\mu_1>0$. 

(i) If $m_1\ge 4$, there is a minimizing solution
of the boundary value problem \re{sus}, \re{sbc} if and only if $\mu_1>m_1-3$.
If $m_1\in\{2,3\}$, there is always such a solution.

(ii) If $m_1\ge 4$ and $\mu_1>(m_1-3)^2/4$, there are even countably many 
$\a:[-\pi/2,\pi/2]\to[-1,1]$ solving \re{sus}, \re{sbc} and
representing smooth Yang-Mills connections on
$\Phi^*TN$, none of which are gauge equivalent.
\end{theorem}

\proof
(i) Minimizing among all $\a$ with $\a(-t)=-\a(t)$, we find a minimizing solution
of the boundary value problem \re{sus}, \re{sbc} if and only if
the constant solution
$\a\equiv 0$ is unstable or $J(0)=\infty$; that is (cf.\
\cite[section~9]{ER1}) if
$\mu_1>m_1-3$ or $m_1\in\{2,3\}$.

(ii) The proof is a very minor modification of the proof of the theorem in
\cite{BC}, where the same is proved for harmonic suspensions $S^n\to S^n$ 
($3\le n\le 6$) of the identity.\qed

\section{Examples of Yang-Mills joins}

Let us see what we can get out of the existence theorems.

\strut\\
{\bf Example 0.} The Levi-Civita connection of $S^{m_1+m_2+1}$ is trivially
Yang-Mills and is the special case $\a(t)=\sin t$ and $\b(t)=\cos t$ (if
$\b$ is needed) that Theorem~\ref{main}, Theorem~\ref{main2}, and 
Theorem~\ref{mainsus} allow if $h_1$ and $h_2$ are identities.

\strut\\
{\bf Example 1.} Nevertheless, Theorem~\ref{mainsus} also produces nontrivial
solutions when applied to $h=\id_{m_1}$ with $4\le m_1\le 8$. And it is
easy to prove that solutions of the o.d.e.\ with different functions
$\a$ cannot be gauge equivalent to each other. This means
that the theorem implies: {\em On every $TS^m$ for $5\le m\le 9$, there are
countably many Yang-Mills connections that are mutually not gauge equivalent.}

\strut\\
{\bf Example 2.} Now we try to join $h_{m_1,\ell}$ ($m_1\ge 2$) with 
$\id_{m_2}$. It is easily checked that the conditions of the existence
theorems are satisfied if $0\le m_2\le 8$. We can formulate that in the 
following way: {\em Each of the pulled back Levi-Civta connections on
$h_{m,\ell}^*TS^{\frac{(2\ell+m-1)(\ell+m-2)!}{\ell!(m-1)!}-1}$ can be 
suspended as Yang-Mills connections $9$ times. }This corresponds to
Smith's observation \cite{Sm} that every harmonic eigenmap can be suspended
harmonically $6$ times.

\strut\\
{\bf Example 3.} The same applies for the case $m_1=1$, that is every
$d_\ell:S^1\to S^1$ can be suspended as a Yang-Mills connection 9 times
(even for $\ell\in\Z$, once we know this for $\ell\in\N$). This is 
geometrically interesting for the following reason: We can interpret the
joined bundles as $f_\ell^*TS^{n}$ for every $\ell\in\Z$ and every 
$n\in\{2,\ldots,10\}$, with $f_\ell$ being a map $S^n\to S^n$ of Brouwer 
degree $\ell$. Depending on $n$, these may be many bundles, maybe
even all $SO(n)$-bundles over $S^n$.

To be more precise, the $SO(n)$-bundles over $S^n$ are classified rather 
easily. Since $SO(n)$ can be covered by just two coordinate patches overlapping
on an annular region around the equator, they are classified by the homotopy
class of the one transition function that is used to patch the two trivial
bundles together; clearly this homotopy class can be seen as an element of
$\pi_{n-1}(SO(n))$. But which element of $\pi_{n-1}(SO(n))$ corresponds to
the bundles $f_\ell^*TS^n$ mentioned above? Since $\pi_{n-1}(SO(n))$ depends
on $n$ in a seemingly unpredictable way, there may be no simple answer. We
can, however, make use of the fact that there is a homomorphism
$e:\pi_{n-1}(SO(n))\to\pi_{n-1}(S^{n-1})\cong\Z$, and that the latter is simply
parametrized by Brouwer degree, hence well-understood. The homeomorphism $e$
is induced by simply evaluating every matrix-valued $A:S^{n-1}\to SO(n)$ 
at some fixed vector $x_0\in S^n$ to give a mapping $Ax_0:S^{n-1}\to S^{n-1}$.

Our first step now is to calculate which element of $\pi_{n-1}(S^{n-1})$ here
corresponds to the tangent bundle $TS^n$. To this end, we observe that
$TS^n$ is parametrized by two coordinate patches 
$f_{\pm}:\overline{B^n}\times\R^n\to TS^n$ given by
\eas
  f_+(x,v)&:=&\Big(x+\sqrt{1-|x|^2}\,e_{n+1}\,,\,
    v-\frac{v\cdot x}{|x|^2}\,x+\frac{v\cdot x}{|x|^2}\Big(\sqrt{1-|x|^2}\,
    \frac{x}{|x|}-|x|e_{n+1}\Big)\Big),\\
  f_-(x,v)&:=&\Big(\xq-\sqrt{1-|x|^2}\,e_{n+1}\,,\,
    v-\frac{v\cdot\xq}{|x|^2}\,\xq+\frac{v\cdot\xq}{|x|^2}\Big(\sqrt{1-|x|^2}\,
    \frac{\xq}{|x|}+|x|e_{n+1}\Big)\Big),
\eeas
where $\xq=(-x_1,x_2,x_3,\ldots,x_n)$. The images $f_+(S^{n-1})$ and 
$f_-(S^{n-1})$ overlap and both parametrize the bundle restricted to the
equator of $S^n$. On $S^{n-1}$, $f_+$ and $f_-$ simplify an read
\eas
  f_+(x,v)&=&(x,v-(v\cdot x)x-(v\cdot x)e_{n+1}),\\
  f_-(x,v)&=&(\xq,v-(v\cdot\xq)\xq+(v\cdot\xq)e_{n+1}).
\eeas
To get $TS^n$, we must find one transition map, and we observe
\[
  f_+(x,v)=f_-(\xq,v-2(v\cdot x)x)=:f_-(\xq,\Phi(\xq)(v)),
\]
where
\[
  \Phi:S^{n-1}\to SO(n),\qquad\Phi(x):=\id-2\xq\otimes\xq
\]
defines the transition map we have been looking for. Now
\[
  \Phi(x)(e_1)=e_1+2x_1\xq,
\]
and it is easily calculated that this mapping $S^{n-1}\to S^{n-1}$ represents
$\pm2\in\Z\cong\pi_{n-1}(S^{n-1})$ if $n$ is even, and $0$ if $n$ is odd.
Similarly, every $f^*TS^{n-1}$ for continuous $f:S^n\to S^n$ can be assigned
an element in $\pi_{n-1}(S^{n-1})$ this way, and this gives a homomorphism
$\pi_n(S^n)\to \pi_{n-1}(S^{n-1})$. Therefore $f_\ell^*TS^n$ represents
a bundle classified by an element of $\pi_{n-1}(SO(n))$ that
$e$ maps to $\pm2\ell\in\pi_{n-1}(S^{n-1})$ if $n$ is even, and $0$ if $n$ is 
odd. Hence we restrict to even $n$ in our search for topologically nontrivial
examples of our Yang-Mills join construction.

Recall that for $n\in\{2,4,6,8,10\}$, we were able to find a Yang-Mills
connection on every $f_\ell^*TS^n$ for every $\ell\in\Z$, and the bundles
correspond to $2\ell\in\Z\cong\pi_{n-1}(S^{n-1})$ under $e$.
The groups
$\pi_{n-1}(SO(n))$ and $\pi_{n-1}(S^{n-1})$ are related by the long exact 
sequence of the homogeneous space $S^{n-1}=SO(n)/SO(n-1)$ which reads
\[
  \ldots\stackrel{e}{\to}\pi_{k+1}(S^{n-1})\to\pi_k(SO(n-1))\to\pi_k(SO(n))
  \stackrel{e}{\to}\pi_k(S^{n-1})\to\pi_{k-1}(SO(n-1))\to\ldots
\]
where $e$ is as before. To illustrate how we can use this, let us first 
consider the case $n=6$. Here is a piece of the exact sequence:
\[
  \begin{array}{ccccccccc}
    \pi_5(SO(5)) &
    \stackrel{0}{\rightarrow} &
    \pi_5(SO(6)) &
    \stackrel{\cdot2}{\hookrightarrow} &
    \pi_5(S^5) &
    \twoheadrightarrow &
    \pi_4(SO(5)) &
    \rightarrow &
    \pi_4(SO(6)). \\
    \Z_2 && \Z && \Z && \Z_2 && 0
  \end{array}
\]
It shows that $\pi_5(SO(6))$ maps injectively to $\pi_5(S^5)$, which means that
every $SO(6)$-bundle over $S^6$ (represented by $j\in\pi_5(SO(6))$) can be
written as $f_j^*TS^6$, and on those we find Yang-Mills connections.
Hence 
\begin{quote}
{\em we have constructed Yang-Mills connections on each of the countably
many (principal) $SO(6)$-bundles over $S^6$. }
\end{quote}
Of course, we always mean $S^6$ equipped with its standard metric. 

The same works for $SO(2)$-bundles over $S^2$, but the result is trivial
because of the sub-critical domain dimension. 

Here are the details for the remaining dimensions. For $n\in\{4,8\}$, we have
\[
  \begin{array}{ccccc}
    \pi_{n-1}(SO(n)) &
    \twoheadrightarrow &
    \pi_{n-1}(S^{n-1}) &
    \rightarrow &
    \pi_{n-2}(SO(n-1)), \\
    \Z^2 && \Z && 0
  \end{array}
\]
hence we find Yang-Mills connections on infinitely many $SO(4)$-bundles
over $S^4$ or $SO(8)$-bundles over $S^8$, but not on all of them. For $n=10$,
\[
  \begin{array}{ccccccccc}
    \pi_9(SO(10)) &
    \stackrel{(\cdot2,0)}{\rightarrow} &
    \pi_9(S^9) &
    \rightarrow &
    \pi_8(SO(9)) &
    \twoheadrightarrow &
    \pi_8(SO(10)) &
    \rightarrow &
    \pi_8(S^9), \\
    \Z\oplus\Z_2 && \Z && \Z_2^2 && \Z_2 && 0
  \end{array}
\]
which shows that we can construct Yang-Mills connections on 
``half of'' the countably many $SO(10)$-bundles over $S^{10}$.

\strut

\strut

\begin{center}
Andreas Gastel\\
Universit\"at Duisburg-Essen\\
Fakult\"at f\"ur Mathematik, Campus Duisburg\\
D-47048 Duisburg, \ Germany\\
{\tt andreas.gastel@uni-due.de}
\end{center}

\end{document}